\documentclass[reqno]{amsart}

\usepackage[english]{babel}
\usepackage[all]{xy}
\usepackage {amssymb}

\usepackage{times}
\usepackage[T1]{fontenc}
\usepackage[latin1]{inputenc}
\usepackage{graphicx}

\theoremstyle{plain}
  \newtheorem{theorem}{Theorem}[section]
  
  \newtheorem*{corollary*}{Corollary}
  
  \newtheorem*{lemma*}{Lemma}
  \newtheorem{proposition}[theorem]{Proposition}
\theoremstyle{definition}
  \newtheorem{definition}[theorem]{Definition}
  
  \newtheorem{observation}[theorem]{Observation}

\DeclareMathOperator{\id}{id}

\DeclareMathOperator{\grad}{grad}

\DeclareMathOperator{\pr}{pr}
\DeclareMathOperator{\ind}{ind}

\begin{document}

\title{
Smooth structure on the moduli space of instantons 
of a generic vector field 
}

\author{Dan Burghelea}

\address{Department of Mathematics,
         The Ohio State University, 
         231 West 18th Avenue,
         Columbus, 
         OH 43210, USA.}

\email{burghele@mps.ohio-state.edu}

\thanks{
    partially supported by NSF grant  no
       MCS 0915996}


\begin{abstract}
This paper is a short version of some joint work with Stefan Haller.  It describes the structure of "smooth manifold 
with corners"  on the space of 
possibly  broken instantons  and completion of unstable manifods of a generic smooth vector field. The result is stated in Theorem 1.4.

 \end{abstract}

\maketitle {}
 
\setcounter{tocdepth}{1}
\tableofcontents

\section{The results} 

This paper has appeared in the volume \cite {Bu5}. Unfortunately we realized that the list of components of the $k-$corner set $\hat {\mathcal M}(x,y)(k)$ of the manifold with corners $\hat {\mathcal M}(x,y)$ considered in Theorem \ref{TH1.5} is incomplete; it misses the components denoted here by  $\hat {\mathcal M}_2(x,y)(k).$  The present version completes this list.
\vskip .1in
 Let $M$ be a smooth closed  manifold and 
 $X: M\to TM$ a  smooth vector field, i.e. a section $X: M\to TM$ in the tangent bundle.
The set  $\mathcal X(X)$ of rest points of $X$  
  consists of points of $M$ where the vector field vanishes, $\mathcal X(X): = \{x\in M  | X(x)=0\}.$   

For any  $x\in \mathcal X(X)$   the differential DX of the smooth map $X$  defines the endomorphism 
$D_x(X): T_x(M)\to T_x (M)$  called  the linearization of $X$ at  $x$
\footnote{$D_x(X)$ is 
defined as follows. Choose an open neighborhood $U$ of $x$ in $M$ and a trivialization of the tangent bundle 
above $U$, $\theta:TU\to U\times T_x(M)$, with $\theta|_{T_x(M)}=\id$. Consider 
$Y:=\pr_2\cdot\theta\cdot X:U\to T_x(M)$ with $\pr_2$ the projection on the second component; 
$Y(x)=0$. Observe that $D_x(Y)$ is independent of $\theta$  and defines $D_x(X):=D_x(Y).$ }.
 
 The rest point $x\in \mathcal X(X)$ is called {\bf  hyperbolic} if the  eigenvalues of $D_x(X) ,$ 
 $\{\lambda \in \text{ Spect} D_x(X)\},$  are complex numbers with real part $ \Re \lambda \ne 0.$  In particular $D_x(X)$ is invertible.

 The hyperbolic rest point $x\in \mathcal X(X)$ is called of {\bf Morse type}   if one can find coordinates  $(u_1, u_2,\cdots, u_n) $ in the neighborhood of $x$ such that  
  $X=\sum _i\pm u_i\frac{\partial}{ \partial u_i}.$
 
Given a hyperbolic rest point $x\in \mathcal X(X)$ the cardinality of the set of eigenvalues  counted with multiplicity 
whose real part is positive  is called {\bf Morse Index} and is denoted by $\ind(x),$
$$\ind (x):= \sharp \{\lambda \in \text {Spect} D_x(X)  |  \Re \lambda >0\}.$$ 
 
 A {\bf trajectory}  of $X$ is a smooth path $\gamma:\mathbb R\to M$ so that $ \frac{d\gamma}{dt}= X(\gamma(t))$ \footnote {If $M$ is not compact then $\mathbb R$ should be replaced by an 
 open interval, the maximal domain of the trajectory; when $M$ is compact this domain is $\mathbb R.$ }.
One  denotes by  $\gamma_y,$  $y\in M,$ the unique trajectory which satisfies  $\gamma_y(0)= y$ and by  $\gamma^-_y$ resp. $\gamma^+_y $ the restriction of $\gamma_y$ to $(-\infty, 0]$ resp. $[0, \infty).$

A trajectory $\gamma$ is called {\bf instanton} from $x$ to $y,$  $x\  \text{and}\ y $ rest points,  if  $\lim _{t\to -\infty} \gamma(t) = x$ and $\lim _{t\to \infty} \gamma(t)  =y.$   
\vskip .1in 

If $x\in \mathcal X(X)$ the set 
$$W^\pm_x:= \{y\in M | \lim_{t\to \pm \infty}\gamma_y=x\}$$ is called the stable / unstable set of $x.$ 

Note that any point $x\in M$ lies  on a trajectory but not necessary on an instanton. For $x,y \in \mathcal X(X)$ 
the set of points lying on instantons from $x$ to $y$    is denoted by $\mathcal M(x,y)$ and is exactly the intersection $$\mathcal M(x,y):= W^-_x\cap W^+_y.$$
The additive group of real numbers $\mathbb R$ acts on $\mathcal M(x,y)$ by "translation" ; if the action is denoted by $\mu:\mathbb R\times \mathcal M(x,y)\to \mathcal M(x,y),$  $a\in \mathbb R, \gamma\in \mathcal M(x,y),$ then $\mu(a, \gamma)(t)= \gamma(t+a).$
The quotient set of this action,  
$$\mathcal T(x,y):= \mathcal M(x,y)/ \mathbb R,$$ is actually the set of instantons from $x$ to $y.$ 

A smooth function $f:M\to \mathbb R$ is called {\bf Lyapunov} for $X$ if 
 $X_m(f)<0, $ for any $m\in M \setminus \mathcal X(X) .$  

A smooth closed one-form $\omega\in \Omega^1(M)$ is called {\bf Lyapunov} if   $\omega(X)_m<0$ for any $x\in M\setminus  \mathcal X(X) .$ 

Not any vector field admits Lyapunov closed one-forms and a vector field can have a Lyapunov closed one-forms but not Lyapunov functions. 
If $X$ has  
Lyapunov functions  then  
any trajectory is an instanton \footnote {This is not true if $M$ is not closed.} and there are no closed trajectories. 
While Lyapunov function  might not exist, for any hyperbolic rest point  
$x\in \mathcal X(X)$  
there exist open  neighborhoods $U$ of $x$ and smooth functions $f: U \to \mathbb R \ $ Lyapunov for $X|_U.$ Similarly, for any instanton $\gamma \in \mathcal T(x,y),$ when  $x,y$ are hyperbolic rest points,   there exist open neighborhoods $U$ of $\gamma$ and smooth functions $f:U\to \mathbb R$  Lyapunov for $X|_U.$
\vskip .1in

The first important result  about stable/unstable sets is the following theorem due to Perron and  Hadamard,  cf. \cite {K95}  Theorem 17.4.3, \cite {AB95} or \cite {I80}  Theorem 6.17 . 

\begin{theorem} If $x$ is a hyperbolic rest point then $W^+_x$ resp. $W^-_x$ is the image of a one to one immersion  $\chi^+_x: \mathbb R^{n-\ind x}\to M$  resp. $\chi^-_x: \mathbb R^{\ind x}\to M.$ 
\end{theorem}

In fact one can find  immersions $\chi_x: \mathbb R^{n-\ind x}\times \mathbb R^{\ind x} \to M$ with $\chi_x(0)=x, $
$ \chi_x |_{(\mathbb R^{n-\ind x}\times 0)}= \chi^+_x$  and
$ \chi_x |_{(0\times \mathbb R^{\ind x})}= \chi^-_x.$

If $X$ admits a Lyapunov function then $\chi^+_x$ resp. $\chi^-_x$ is actually  a smooth embedding which makes $W^+_x$ resp. $ W^-_x$ a  smooth submanifold of $M.$ 
In general the topology of $W^+_x$ resp. $ W^-_x$ obtained by  identification with $\mathbb R^{n-\ind x}$ resp. $\mathbb R^{\ind x},$ and referred below as {\it manifold topology}, is not necessary
the same as the induced topology \footnote{ The later being  coarser in general.}. 
The immersion $\chi^\pm_x$ is not unique but the manifold topology on $W^\pm_x$ and the smooth structure defined by the chart $\chi^\pm_x$  is. It is possible  that for a hyperbolic rest point $x$ the stable set and the unstable set  be the same  
\footnote {In this case $\dim M$ has to be even.},  but  the manifold topologies are  different. 
 
In conclusion   $W^\pm_x$ is  equipped with a canonical structure of smooth manifold  such that the canonical inclusion $i^\pm_x: W^\pm_x\to M$ is a one to one immersion (embedding if $X$ admits Lyapunov functions).  
Suppose $x,y\in \mathcal X(X)$ are hyperbolic and the maps $i^-_x$ and $i^+_y$ are transversal.
Then the set 
$$\mathcal M(x,y):= \{(u,v)\in W^-_x\times W^+_y | i_x^-(u)= i_y^+(v)\}$$
has a structure of a smooth manifold of dimension $\ind x- \ind y$
with the canonical inclusion  $i_{x,y}:\mathcal M(x,y) \to M$ a one to one smooth immersion and  the action $\mu:\mathbb R\times \mathcal M(x,y)\to \mathcal M(x,y)$  smooth and free  provided $x\ne y$. In this case the  quotient space $\mathcal T(x,y)$  receives a canonical structure of smooth manifold  of dimension $\ind x-\ind y -1$ with the quotient map $p:\mathcal M(x,y)\to \mathcal T(x,y)$ a smooth bundle.  
Clearly if  $\ind x\leq \ind y$ and $x\ne y$ then $\mathcal M(x,y)$ is empty.

\vskip .1in
From now on we  suppose  the vector fields satisfy  the following two properties :

{\bf $P_1:$} All rest points of $X$ are hyperbolic.

{\bf $P_2:$} For any two rest points $x, y\in \mathcal X(X)$ the maps   $i^-_x$ and $i^+_y$ are transversal. 
\vskip .1in 
The following result due to Kupka and Smale,  cf \cite {K63}, \cite {S63}, \cite {P67}, shows that this  is the  generic situation.
\begin{theorem} \label{th0}
The set of vector fields which satisfy $P_1$ and $P_2$ are residual \footnote { i.e. contain a countable intersection of open dense sets}  in the $C^r-$ topology for any $r\geq 1.$
\end{theorem} 
\hskip 2in $\therefore$
\vskip .1in

Write  $x>y$ for  $\ind x > \ind y$ when  $x,y\in \mathcal X(X).$     
\vskip .2in
For any $k\geq 2,$  introduce
$$
\hat W ^-_x (k):= \bigsqcup _{\tiny{
\begin{array}{cccc} 
x_1>x_2>\cdots >x_k \mid
 x>x_1
 \end{array}}}  
\mathcal T(x,x_1) \times \mathcal T(x_1, x_2)\times (\cdots  \mathcal T(x_{k-1}, x_k) \times W^-_{x_k},$$
\vskip .1in
$$
\hat W ^+_y (k):= \bigsqcup _{\tiny{
\begin{array}{cccc} 
y_k> y_{k-1}>\cdots >y_1 \mid
 y_1>y
 \end{array}}}  
\mathcal  W^+_{y_k}\times T(y_k, y_{k-1}) \times \cdots \mathcal T(y_2, y_1)\times  \mathcal T(y_1, y) ,$$
\vskip .1in
$${\hat{\mathcal T}(x ,y)(k):= \bigsqcup _{\tiny{\begin{array} {cccc} x_1> x_2>\cdots > x_{k-1}\mid
 x>x_1, x_{k-1}>y\end{array}}} \mathcal T(x_0,x_1)\times \mathcal T(x_1,x_2) \cdots \mathcal T(x_{k-1}, y)},$$
\vskip .1in
$$\hat {\mathcal M}_1(x ,y)(k):=
 \underset{ \tiny{\begin{array} {cccc} r=0,1,2,\cdots, k\\
 x_0>x_1>\cdots >x_{k+1} \\
\rm{with}\ x=x_0, y= x_{k+1} \end{array}}} 
{\bigsqcup} 
\begin{cases}&\ \mathcal T(x_0, x_1)\times \cdots  \mathcal T(x_{r-1},x_r)\times \mathcal M(x_r , x_{r+1})\times 
\\ &\times \mathcal T(x_{r+1},x_{r+2})\times \cdots \times \mathcal T(x_k,x_{k+1})
\end{cases}$$
when $r=0$ resp.  $r=k$ the formula begin with $\mathcal M(x_0,x_1)$ resp. ends up with $\mathcal M(x_k, x_{k+1})$
\vskip .2in
$$\hat{\mathcal M}_2(x ,y)(k):= 
 \underset{ \tiny{\begin{array} {cccc}r=0,1,2,\cdots, k\\
x_0> x_1>\cdots >x_{k} \\
\rm{with}\ x=x_0, y=x_k \end{array}}} 
{\bigsqcup} \begin{cases}&\mathcal T(x_0, x_1)\times \cdots  \mathcal T(x_{r-1},x_r)\times \mathcal M(x_r , x_{r})\times \\
&\times \mathcal T(x_{r},x_{r+1})
\cdots \times \mathcal T(x_{k-1},x_k)
\end{cases}$$
\newline when $r=0$ resp. $r=k$ the formula begins with $\mathcal M(x_0,x_0)$ resp. ends with $\mathcal M(x_k, x_k)$
\vskip .2in
For $k=1$
$$\hat W^-_x(1)= \bigsqcup _{\tiny{ x_1\mid x>x_1}} \mathcal T(x,x_1)\times W^-_{x_1}$$
$$\hat W^+_y(1)= \bigsqcup _{\tiny{ y_1\mid y_1 >y}} W^+_{y_1}\times \mathcal T(y_1,y)$$
$$\hat {\mathcal T}(x,y)(1)= \bigsqcup _{\tiny{ y_1 \mid x> y_1 >y}} \mathcal T(x,y_1)\times \mathcal T(y_1,y)$$
$$\begin{aligned}\hat {\mathcal M}_1(x,y)(1)= &\bigsqcup _{\tiny{ x_1 \mid  x>y_1 >y}} \mathcal M(x,x_1)\times \mathcal T(x_1,y) \sqcup \bigsqcup _{\tiny{ x_1 \mid  x>x_1 >y}} \mathcal T(x,x_1)\times \mathcal M(x_1,y)\\
\hat {\mathcal M}_2(x,y)(1)= &\quad \quad \mathcal M(x,x)\times \mathcal T(x,y) \quad \sqcup \quad \mathcal T(x,y)\times \mathcal M(y,y)
\end{aligned}$$
\vskip .2in
For $k=0$
$$\hat W^-_x(0)=W^-_x, \hat W^+_y(0)=W^+_y,\ \   \hat {\mathcal T} (x,y)(0)= \mathcal T (x,y), \ \  \hat{\mathcal M} (x,y)(0)= \mathcal M (x,y).$$
\vskip .2in
For $k>0$ define the maps 
$\hat i ^{\mp}_x(k): \hat W ^{\mp}_x(k)\to M$ and  $\hat i_{x,y}(k): \hat{\mathcal M}(x,y)(k)\to M$ by 
$$ \hat i ^{\mp}_x(k)= \bigsqcup _{\tiny{
\begin{array}{cccc} 
y_0< y_1<\cdots <y_k \\
 x=y_0
 \end{array}}}  \hat i ^{\mp}_{y_k}\circ pr_{W^{\mp}_{y_k}}, $$

 $$
 \hat {i _1} _{x,y}(k)= \underset{  \tiny{\begin{array} {cccc}y_0> y_1>\cdots >y_{k+1}
\\
 x=y_0, y_{k+1}=y \end{array}}} {\bigsqcup}  i_{y_i,y_{i+1}}\circ pr_{\mathcal M_{y_i,y_{i+1}}}.$$

 $$
 \hat {i _2}  _{x,y}(k)= \underset{  \tiny{\begin{array} {cccc}  r= 0,1,\cdots ,k\\
y_0> y_1>\cdots y_{r} > \cdots y_{k}\\
x=y_0, y_{k}=y \end{array}}} {\bigsqcup} i_{y_r}\circ pr_{\mathcal M_{y_r,y_r}}.$$
 
 Here $i_y$ denotes the inclusion $y\in M$

and for $k=0$ 

 $$\hat {i }^{\mp}_x(0)= i_{x}, \ \ \hat i_{x,y}(0)=  i_{x,y}.$$

Denote by  $\hat {W}^{\mp}_x, \ \hat {\mathcal T}(x,y)$ and $\hat {\mathcal M}(x,y)$ the sets defined by
$$\hat {W}^{\mp}_x:= \bigsqcup_{k\geq 0} W^{\mp}_x(k), $$ 
$$\hat {\mathcal T}(x,y):= \bigsqcup_{k\geq 0} \mathcal T(x,y)(k),$$  
$$\hat {\mathcal M}(x,y):= \bigsqcup_{k\geq 0} \hat{\mathcal M}_1(x,y)(k) \sqcup \bigsqcup_{k\geq 0}\hat{ \mathcal M}_2(x,y)(k) ,$$  
and  denote by $\hat i^{\mp}_x: \hat W^-_x\to M$ and $\hat i_{x, y}: \hat {\mathcal M}(x,y)\to M$ the maps defined by
$$\hat i^{\mp}_x |_{\hat W ^{\mp}_x(k)} := \hat i ^{\mp}_x(k),$$
$$ \hat i_{x,y}|_{\hat{\mathcal M}(x,y)(k)} :=  \hat i _{x,y}(k).$$

We equip  $\hat W^\mp_x, \hat {\mathcal T}(x,y)$ and $\hat{\mathcal M}(x,y)$ with the transversal slice topology defined below. 
Since $\hat W^-_x$ for the vector field $X$ is the same as $\hat W^+_x$ for the vector field $-X$ and since 
 $$ \hat{\mathcal M}(x,y)= \{(u,v)\in \hat W^-_x\times \hat W^+_y |  (\hat i^{-}_x)(u)= \hat  i^+_y(v)\},$$
it suffices to describe only the transversal slice topology for  $\hat W^-_x$ and $\mathcal T(x,y).$   For this purpose a few more definitions are  necessary.  

A   {\bf broken instanton} $\underline {\gamma}= (\gamma^1, \gamma^2, \cdots \gamma^k)$ consists of a collection of instantons ${\gamma^i}$     's  with the property 
$$\lim _{t\to \infty}\gamma^i(t)= \lim _{t\to -\infty}\gamma^{i+1}(t), \ \ i=1,2\cdots, k-1.$$ 

An element of  $\hat{\mathcal T}(x,y)$ is a broken instanton with the property 
$$\lim _{t \to- \infty}\gamma^1(t)=x,  \lim _{t \to \infty}\gamma^{k}(t)=y.$$

An element of $\hat W^-_x$   can be uniquely represented  as a pair $\tilde\gamma :=(\underline \gamma,  m),$ with $\underline \gamma= (\gamma^1, \gamma^2, \cdots \gamma^k)$ 
a broken instanton and  $m\in M$ which satisfy:

\hskip .1in 1. $\lim _{t \to- \infty}\gamma^1(t)=x, $ 

\hskip .1in 2. $\lim _{t \to \infty}\gamma^{k}(t)=\lim_{t\to -\infty}\gamma_m(t).$

Given a collection of trajectories  $\underline{\gamma}
= (\gamma^1, \gamma^2,\cdots \gamma^k)$  
with the property 
$\lim _{t\to \infty}\gamma^i(t)= \lim _{t\to -\infty}\gamma^{i+1}(t) \ ,$
i.e. a { broken trajectory}, a {\bf  transversal slice} is
a collection $\underline U= (U_1, U_2,\cdots U_r)$ of disjoint $(n-1)$ dimensional submanifolds diffeomorphic  to open  discs  
which are  {\it transversal}  to $X$ (transversal to trajectories of $X$)  and satisfy:

\hskip .1in 1. each $\gamma_i$ intersects at least one $U_j,$ 

\begin{figure}[ht] \label{F;1}
\begin{center}
\includegraphics[height=60mm]{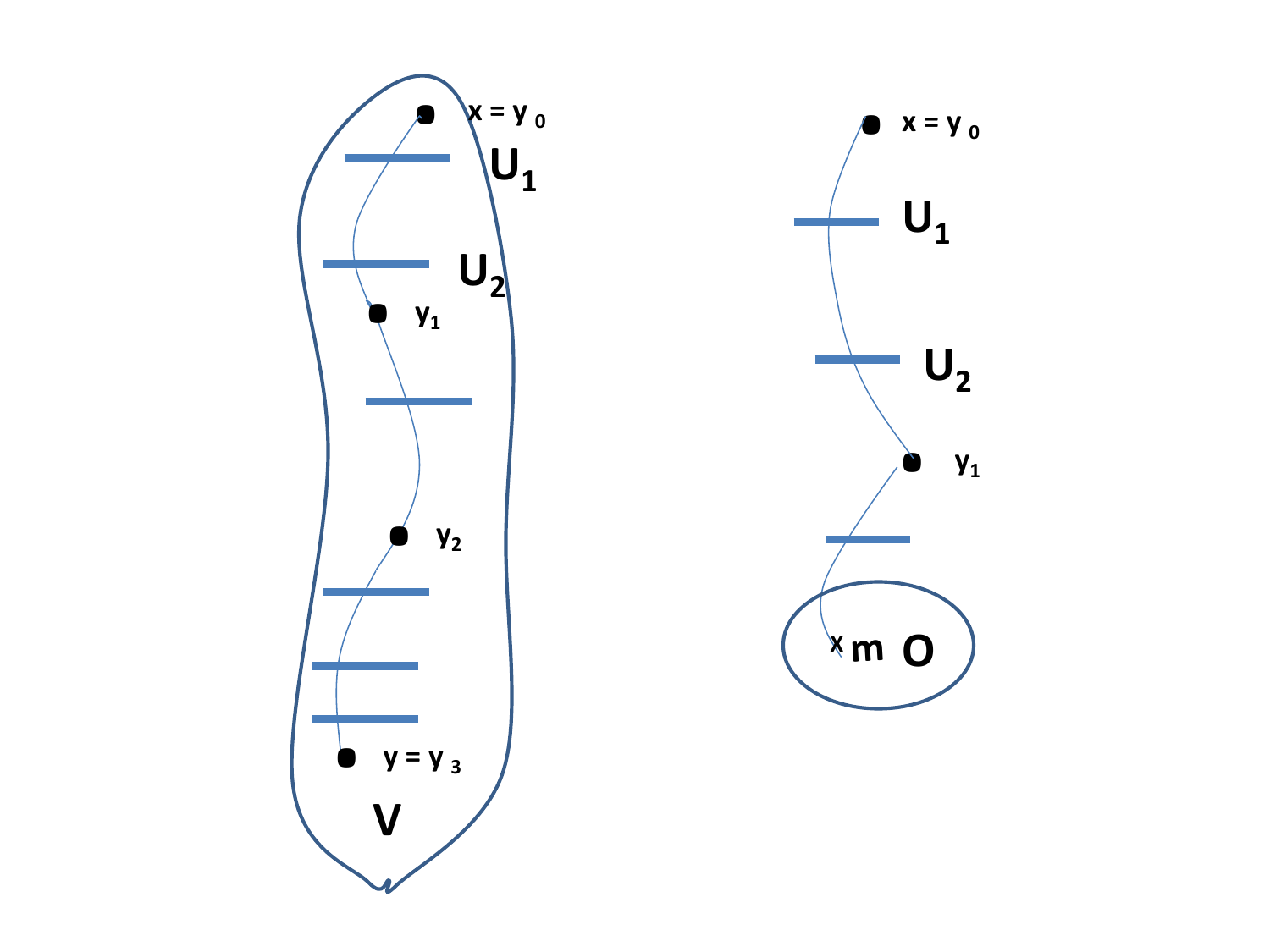}
\caption{}
\label{fig:DANU7}
\end{center}
\end{figure}

\hskip .1in 2. consecutive $U_j'$s  intersect either the same  or consecutive $\gamma^i$'s, hence $ r\geq k,$   
 cf. FIGURE \ref{F;1}. 
\vskip .1in

For a system $(\underline \gamma, V, \underline U)$  with $V$ open neighborhood of $\underline \gamma$ in $M$ 
and  $\underline U$  a transversal slice to $\underline \gamma\ $  denote by $\mathcal V(\underline \gamma, V, \underline U)$ the collection of  
possibly broken instantons  in $\hat{\mathcal T}(x,y)$  which lie inside $V$ and have $\underline U$ as a transversal slice.

For a  system  $(\tilde \gamma, V,  \underline U, O)$  with $\tilde\gamma\in  \hat W^-_x,$ $V$ open neighborhood of  $\tilde \gamma$ in $M,$  
$ \underline U$ transversal slice for  $(\gamma^1,\cdots \gamma^r, \gamma^-_m)$ 
and $O$ an open neighborhood of $m,$ see FIGURE \ref{F;1},
denote by $\mathcal V(\tilde \gamma, V, \underline U, O)$ 
the collection of all elements in $  \hat W^-_x$  which lie inside $V,$ have $\underline U$ as a transversal slice  and have the end point in $O.$ 

A base of the {\bf transversal slice topology} for $\hat {\mathcal T}(x,y)$ resp. for $\hat W^-_x$ 
is provided by the sets $\mathcal V(\underline \gamma, V,  \underline U)$ for all systems  $(\underline \gamma, V, \underline U)$ resp. 
the sets $\mathcal V(\tilde \gamma, V, \underline U, O)$  for all systems  $(\tilde \gamma,\  V,  \underline U, \ O).$

\begin{theorem}\label {TH1.3}

1. When equipped  with the "transversal slice"  topology  the sets $\hat W^\mp_x, \hat { \mathcal T}(x,y)$ and  $ \hat { \mathcal M}(x,z)$
are Haussdorf paracompact spaces and the maps $\hat i^-_x: \hat W^-_x\to M$ and  $\hat i_{x,z}: \hat {\mathcal M}(x,z)\to M$ 
are continuous.

2. If $X$ admits Lyapunov function  then $\hat W^\mp_x, \hat { \mathcal T}(x,y)$ and $\hat{\mathcal M}(x,z)$
are compact. 
\end{theorem} 

\proof (sketch).  Suppose first  that $X$ admits a Lyapunov function $f:M\to \mathbb R.$ 

{\it The case of $\hat{\mathcal T}(x,y)$:}  Suppose $f(x)= c_1$ and $f(y)=c_2$ and choose $\alpha_0<c_1 <\alpha_1 <  \alpha_2<\cdots \alpha_{k-1}< c_2 < \alpha_k$ with $\alpha_i$ regular values and with the  intervals $(\alpha_i, \alpha_{i+1}) $ containing only one  critical value $c_{i+1}.$  The map 
which assigns to a broken instanton $\underline \gamma$ its intersection with the levels $f^{-1}(\alpha_i)'s,$
$$i:\hat{\mathcal T}(x,y)\to \prod _{i= 1,\cdots, (k-1)} f^{-1}(\alpha_i),$$ 
 is a one to one map.  It is not hard to see that the transversal slice topology is the same as the  topology induced by this embedding. Indeed one can consider only transversals which lie in the levels $f^{-1}(\alpha_i)$ and show they suffice to describe the transversal slice topology.  Standard arguments cf. \cite{BH01} show that the image of $i$ is closed hence compact.  This proves the result for  $\hat{\mathcal T}(x,y).$ 

{\it The case of $\hat W^-_x: $} 
Write the critical values in decreasing order  $\cdots > c_i>c_{i+1}>\cdots .$  Verify the first assertion in Theorem \ref {TH1.3} for   $(\hat i^-_x)^{-1}(f^{-1} (a, b))$ instead of $\hat W^-_x$ when  $a, b$ are regular values of $f$ with at most one critical value in the interval $(a,b).$  The set $(\hat i^-_x)^{-1}(f^{-1} (a, b))$
can be embedded  in a product of  finitely many levels of   $f$ and  $M$ and one can  check that the transversal slice topology and the  topology induced by such  embedding  are the same. 
Since  $\hat i^-_x$ is continuous w.r. to the transversal slice topology, hence $(\hat i^-_x)^{-1}(f^{-1} (a, b))$ is open and $\hat W^-_x$ is covered by such sets, then  the conclusion extends from $(\hat i^-_x)^{-1}(f^{-1} (a, b))$ to $\hat W^-_x.$  The compacity assertion follows   from the compacity of $(\hat i^-_x)^{-1}(f^{-1} [\alpha, \beta]).$

 To conclude the statement in  general  (when no Lyapunov function exists) one observes  that  any  $\underline \gamma \in \hat {\mathcal T}(x,y),$   or  $(\tilde \gamma= ( \gamma^1,\cdots \gamma^r, \gamma^-_m) \in \hat W^-_x$  lies inside an open set $V$ of $M$ so that the vector field $X|_V$ admits Lyapunov function $f:V\to \mathbb R.$ This follows  from the existence of Lyapunov functions in the neighborhood of each hyperbolic rest point a fact noticed above.  
More details will be  contained in  \cite {BH03}.  

The main result of this paper states that the topological spaces $\hat W^\pm_x, \  \hat {\mathcal T}(x,y)$ and $\hat{\mathcal M}(x,y)$  have structures of smooth manifold with corners. To explain this let us recall a few definitions. 
\vskip .1in

The standard example and  the local model of a smooth manifold with corners is 
$$\mathbb R^n_+ :=\{(x_1, x_2, \cdots x_n) \in \mathbb R^n | x_i \geq 0\}.$$

The $k-$corner of $\mathbb R^n_+$ is 
$$ {\partial_k \mathbb R^n_+:= \{x\in  \mathbb R^n_+ | \text {exactly} \ k \  \text{coordinates} \ x_i=0\} }.$$

Denote by $[[n]]$ the set $\{1,2,\cdots, n-1, n\}.$ If   $I=\{ {i_1},
{i_2},\cdots {i_k}\}$ is a subset of $[[n]]$   denote  by $\mathbb R^I$
the set of points in $\mathbb R^n_+$ whose coordinates $x_{i_1},
x_{i_2},\cdots x_{i_k}$ are different from zero while all other coordinates
vanish.  Note that each $\mathbb R^I$ carries a canonical orientation
defined by the order $i_1 <i_2,\cdots< i_k.$   

 A {\bf smooth manifold with corners} is   a Haussdorf paracompact space  $P$ 
equipped with a {\it differential structure}
locally isomorphic to $\mathbb R^n_+.$ A differential structure is given by an equivalence class of atlases.   
An atlas $\{\varphi_\alpha:  U_\alpha \to V_\alpha, \alpha \in \mathcal A\}$ consists of 

an index  set $ \mathcal A$ ,  

$U_\alpha \subseteq P$ open sets of $P,$   
 
 $V_\alpha \subseteq \mathbb R^n_+$  open sets of $\mathbb R^n_+$,  and
 
 $\varphi_\alpha $ homeomorphisms (charts)
 
 \noindent so that $\cup_\alpha  V_\alpha= P$  
 and $\varphi_\beta\cdot \varphi^{-1}_\alpha$  
 are smooth and of maximal rank where defined. Two atlases are equivalent if their data  considered together remain an atlas.
 
 The $k-$corner $\partial_kP$ is the set of points which in some chart (and then in any)  
correspond to $\partial_k \mathbb R^n_+.$

The manifold with corners is orientable if $\partial _0P$ is orientable.  
An orientation on
such manifold is an orientation for  its  tangent bundle, equivalently an orientation of the open manifold $\partial_0P.$ 
. 

A smooth manifold with corners $P$ is {\bf clean}  if the closure of each connected component of the corners  is a smooth manifold with corners. 
\vskip .1in
The main result of this paper is the following: 

\begin{theorem}\label{TH1.5}

Let $X$ be a smooth vector field satisfying $P_1$ and $P_2$ (defined before Theorem \ref{th0}).

1.  There exists  a canonical  structure of clean smooth manifold 
with corners on  
$ \hat W^\mp_x, \ $ $ \hat { \mathcal T}(x,y)$ and $ \ \hat { \mathcal M}(x,y)$  with  
$$ \partial_k \hat W^-_x = \hat W ^-_x(k), \ \ 
\partial_k \hat { \mathcal T}(x,y) = \hat{\mathcal T}(x,y)(k), \ \ 
 \partial _k \hat { \mathcal M}(x,y)= \hat{\mathcal M}(x,y)(k)$$ and $\ \hat i^\mp_x$  so that  $\ \hat i_{x,y}$ are smooth maps
 \footnote {A previous version of this paper  was published in \cite{Bu3}. Unfortunately the description the $k-$corner of $\partial _k\hat{\mathcal M}(x,y)$ was incomplete. It misses the term $\hat{\mathcal M}_2(x,y)(k)$ which was added in this version. apparently no other additions /corrections were necessary in the published text}. 

2. If  the rest points of $X$ are of Morse type 
there exists an additional structure of smooth manifold with corners
(different but  diffeomorphic to the structure stated in 1.) with the same
corners  and the identity map restricted to each $k-$corner a  diffeomorphism.

3. Both $\hat{\mathcal T}$ and  $\hat{\mathcal M}$ are equipped with (stable)
framings.  A collection of orientations $\{ o_x, x\in \mathcal X(X)\}$ on $W^-_x$
induces coherent  orientations on $\hat{\mathcal M}(x,y)$ and $\hat{\mathcal T}(x,y)$
\footnote{ This means that for any three rest points $x,y,z$ with $\ind x >\ind y >\ind z$ the orientation $o_{x,z}$ on $\mathcal T(x,z)$ induces on $\mathcal T(x,y)\times \mathcal T(y,z) \subset \partial \hat{\mathcal T}(x,z)$ the oposite of the orientation $o_{x,y} \otimes o_{y,z}.$}.
\end{theorem}
 
Theorem \ref{TH1.5} is not new  but in the generality formulated above  inexistent in literature.  In less generality it can be recovered from \cite{La} and \cite{AB95}  for the gradient of a  Morse function and from  \cite {BFK}, \cite {BH01}, \cite{BH02} or \cite{La} for the gradient of a closed one form. The proof below  is along  the lines of \cite{BFK} or \cite{BH01}.
In  \cite{BH03} the result will be proven  for a more general class of vector fields, called  HB  (hyperbolic - Bott) vector fields. For these vector fields  the  set of rest points $\mathcal X(X)$  is  a  smooth submanifold with $D_x$ hyperbolic in normal directions of  $\mathcal X(X).$

The  smooth structure  provided by 
Theorem \ref{TH1.5}  is not the only possible canonical smooth structure.   
In fact, in the case the  rest points are of  Morse type, Theorem \ref{TH1.5} provides two such canonical structures never the same. However all  smooth structures  of manifold with corners on $\hat W^\mp_x, \ \hat{\mathcal M} (x,y)$ or $
 \hat{\mathcal T}(x,y)$ which extend the smooth structure of $  \hat W ^\mp_x, \ \hat {\mathcal M} (x,y)$ or $ \hat {\mathcal T}(x,y)$  are  diffeomorphic.  By elementary smoothing theory one can show that  such diffeomorphisms  can be chosen to be the identity on arbitrary closed subsets  of the  $\partial_0 $ part. 
\vskip .1in

Theorem \ref{TH1.5} provides a source of new invariants which deserve
attention and we plan to explore in future work. For example:

1. A Morse type complex can be assigned to a  class of
vector fields substantially larger  than the gradient like vector fields. Its
homology/cohomology referred to as the {\it instanton homology
 (cohomology)} might relate the topology of the manifold and the dynamics of $X$ in
a  more subtle way than in the case of gradient like vector fields. For example if $X$ is the gradient of a closed one form  both
the Novikov  cohomology 
and the cohomology of $M$ twisted by a closed one form  can be
obtained as  instanton homology/cohomology.  More general vector fields lead to more subtle instanton homologies / cohomologies.

2. A  chain/ cochain complex can be  derived from the corners structure of
the manifold with corners $\hat{\mathcal T}(x,y).$  The homology/cohomology of such a complex, referred to as the
 {\it incidence cohomology} seems  natural to investigate. It carries significant dynamical information not obviously related 
 to the topology of the manifold.

3. The stable framing of $\hat{\mathcal T}(x,y)$ can be used to define 
elements in $ \pi_i^S(\Omega M),$ the  stable homotopy groups of the
free loop space of $M.$  A parametrized version of such elements might provide a more analytic  understanding of the relationship between the homotopy of the space of diffeomorphisms and the Waldhaussen K- theory of the underlying manifold.

{\it Acknowledgement:}  I thank Stefan Haller for pointing out a number of errors and misprints in a previous version of this manuscript.

\section{Some basic ODE}
Recall  that a linear transformation of $\mathbb R^n$  is hyperbolic if all its eigenvalues  have non vanishing real part.  The stable resp. unstable subspace, $\mathbb R^+\subseteq \mathbb R^n$ resp. $\mathbb R^-\subseteq \mathbb R^n ,$ are the sum of generalized eigenspaces corresponding to the eigenvalues with negative resp. positive real part.

Consider
\begin{enumerate}
\item $A\in M(n\times n; \mathbb R)$  hyperbolic  with stable space 
 $\mathbb R^+ = \mathbb R^k\times 0 $ and unstable space  
$\mathbb R^-= 0\times \mathbb R^{n-k},$ 
\item $g: \mathbb R^n \to  \mathbb R^n= \mathbb R^+\times \mathbb R^-$
a smooth map  with compact support
with $g(0)= 0$ and $ D_0g=0.$ We write $g= (g^+, g^-).$ We require in addition that $g^+(x^+, 0)=0, g^-(0, x^-)=0.$
\end{enumerate}

Let $X:  \mathbb R^n\to \mathbb R^n$ be the  smooth map 
defined by  
$X(x)= Ax+g(x).$   
Regard $X$ as a smooth vector field on $\mathbb R^n.$ Clearly $0\in \mathbb R^n$ is a hyperbolic rest point and the only rest point  in a small neighborhood of $0.$ 

Denote by
 $$\gamma (t: p, q, T_1, T_2):= (\gamma^+(t: p, q, T_1, T_2), \gamma^-(t: p, q, T_1, T_2))\in \mathbb R^+\times \mathbb R^-$$
 $(p, q)\in \mathbb R^+\times \mathbb R^-,$  $T_1 < T_2$ 
a  trajectory which satisfies  
\vskip .1in
\begin{equation}
\begin{aligned}
\gamma^+(T_1: p, q, T_1, T_2)= p,\\
\gamma^-(T_2: p, q, T_1, T_2)= q.
\end{aligned}
\end{equation}

In general such  trajectory might not exist and even if exists it might not be unique.
However we have :

\begin{theorem}  \label{T;2.1} For any positive integer  $N$ there exists  $\epsilon, \rho,  C>0$ so that :
\vskip .1in
1.  For any $(p, q) \in D^+(\epsilon)\times  D^-(\epsilon),$ the discs of radius $\epsilon$ in $\mathbb R^+\times R^-,$  and any $T_1 < T_2$  
  there exists a unique trajectory $\gamma(t: p, q, T_1, T_2).$
  \vskip .1in

2. Moreover the following estimates hold
\begin{equation}
\begin{aligned}
||\mathbb T\gamma^+(t: p,q,T_1, T_2)|| \leq  \epsilon  C^{-\rho (t-T_1)}\\
||\mathbb T\gamma^-(t: p,q,T_1, T_2)|| \leq  \epsilon  C^{\rho (t-T_2)} 
\end{aligned}
\end{equation}
for $T_1 <t< T_2, $ where  $\mathbb T= \frac{\partial^{k_1} }{\partial p^{k_1}} \frac{\partial ^{k_2}}{\partial q^{k_2}} \frac{\partial ^{k_3}}{\partial t^{k_3}} \frac{\partial ^{k_4}}{\partial T_1^ {k_4}} \frac{\partial^{k_5} }{\partial {T_2}^{k_5}}$
with $k_1+k_2 + k_3 +k_4+k_5 \leq N$ 
\end{theorem}

The result is a straightforward application of contraction principle.   A proof can be derived  on the lines of the proof of  
Theorem A.2 and Lemma A.3 of  Appendix of  \cite{AB95}.  For the reader's convenience we sketch  the  proof  of (1.) and comment on the proof of (2.). 

We continue to write
 $ p$ and $ g^+$ for  $(p,0),$ and  $(g^+,0),$ 
 and $ q$ and $
  g^-$ for  $(0,q)$ and 
$(0, g^-).$ 
We also  write
$\Psi(s,t)$ for $e^{(t-s)A}.$

Using the hyperbolic linear transformation $A$ one can produce the  real numbers $\rho'>0$ and $C >1$ so that 
\begin{equation}
\begin{aligned}
||\Psi(s,t) p|| &\leq Ce^{-\rho' (s-t)}|| p||\  \text{for}\  s\geq t,   \\
||\Psi(s,t)  q|| &\leq Ce^{\rho' (s-t)}|| q||\  \text{for}\  s\leq t.
\end{aligned}
\end{equation}

Since  the trajectory $\gamma (t):= \gamma (t, p,q, T_1, T_2)$ has to satisfy  the equality 
\begin{equation}\label{ES}
\gamma (t)= \Psi(t, T_1)(p) +   \Psi(t, T_2)(q) + \int_{T_1}^t\Psi(s,t)g^+(\gamma (s))ds  - \int_t^{T_2} \Psi(s,t) g^-(\gamma(s))ds 
\end{equation}
one concludes that   $\gamma (t)$ is a fixed point of the map $F: C^0([T_1, T_2], \mathbb R^n) \to C^0([T_1, T_2] ,\mathbb R^n) $ \footnote {
 $C^0([T_1, T_2], \mathbb R^n)$ denotes the Banach space of continuous function from $[T_1, T_2] $ to $\mathbb R^n$ with the $C^0-$ norm.}
defined by 
\begin{equation}\label {E:3}
F(x(t))= \Psi(t, T_1)( p) +   \Psi(t, T_2)( q)  + \int_{T_1}^t\Psi(s,t) g^+(x(s))ds  - \int_t^{T_2} \Psi(s,t)(0,  g^-(x(s))ds.
\end{equation}
Choose $B >0$ so that $$||Dg (x)||\leq B ||x||.$$

Then if $||(p,q)||\leq \epsilon$ and $||x(t)||\leq \eta$
one has 

\begin{equation}\label {E6}
 ||F(x(t)|| \leq 2C \epsilon + \frac{2 CB \eta^2}{\rho'} .
\end{equation}
One can  find  $\eta$ and $\epsilon'$  small so that  $F$ sends the disc of radius $\eta$ into itself
provided $||(p,q)||\leq \epsilon'$.
Precisely  
one 
chooses $\eta$  to satisfy
 $$\eta \leq  \frac{\rho'}{4BC} ,$$
 and $\epsilon'$ to satisfy $$\epsilon'  < \frac{\eta}{4C}.$$ 
These choices make both terms of the right side  of the above inequality (\ref{E6}) smaller than $\eta/2.$
Then  $||x(t)||\leq \eta$ implies  $||F(x(t)) ||\leq \eta.$

Note that the  estimates remain true for any $\eta' \leq \eta$ (when $\epsilon' $ is appropriately chosen).
 Since we have 
 \begin{equation}
||\frac{\partial F}{\partial x} ||   \leq  CB  \eta/2\rho'  (1-e^{\rho'(T_1-t)}) + CB\eta/ 2\rho'  \leq CB  \eta/ \rho' < 1/4
 \end{equation}

 \noindent by choosing $\eta<\rho'/4CB$ and $\epsilon' <\eta/4C$ one concludes that $F$ sends the disc $\mathbb D(\eta)$ of radius $\eta$ (in the complete metric space $C^0([T_1, T_2], \mathbb R^n)$) into itself and that $F: \mathbb D(\eta)\to \mathbb D(\eta)$ is a contraction. 

To prove (2.) we first produce  $\epsilon $, $\rho$  with $\epsilon<\epsilon'$, $\rho<\rho'$   so that inequalities (2) are satisfied for $N=0,$ 
then decrease $\epsilon $ and $\rho$  inductively  
with $N$  to satisfy (2.) 
 
This is done by incorporating  the estimates  stated in (2.) in the definition of the metric space  $C^0([T_1, T_2], \mathbb R^n),$ 
and decrease $\epsilon, \rho$  to make sure that $F$ remains a contraction even in the presence of these estimates, hence the unique fixed point $\gamma $ satisfies  the estimates. 

\section {Elementary differential topology of smooth manifolds with corners}

If  $P$ is  a smooth manifold with corners then:

\hskip .2in $\partial _k P$ is a a smooth $(n-k)-$manifold,

\hskip .2in $\partial P:= \bigsqcup_{k\geq 1}\partial_k P$ is a topological $(n-1)-$manifold, and

\hskip .2in $(P,\partial P)$ is a smoothable topological manifold with boundary.

with an unique smooth structure up to a diffeomorphism

If $P_i,$ $i=1,2$ are two smooth manifolds with corners then the  product $P_1\times P_2$ is a smooth manifold with corners with 
$$\partial_k (P_1\times P_2)=\bigsqcup_{r=k'+k''} \partial _{k'}P_1\times \partial_{k"}P_2.$$ 
If both $P_1$ and $P_2$ are clean then so is the product.

Let $P$ be a smooth  manifold with corners, 
$M$ a smooth manifold,  $S\subseteq M$ a smooth submanifold 
and $f:P\to M$ a  smooth map.

\begin{definition} 
The map $f $ is transversal to $S$,  written $ f\pitchfork  S $,   if $f|_{\partial _k P} \pitchfork  S$ for any $k=0,1,\cdots n.$ 
\end{definition}

Let $P_i,$ $i=1,2$ be two smooth manifolds with corners, $M$ a smooth manifold and $f_i: P_i\to M$ smooth maps.
 \begin{definition}
The  maps $f_i: P_i\to M$  are transversal, written  $f_1\pitchfork f_2,$
if the product $f_1\times f_2 : P_1\times P_2\to M\times M$ is transversal to the diagonal $\Delta_M:=\{(x,x)\in M\times M | x\in M\}$
 $(f_1\times f_2 \pitchfork \Delta_M).$
 \end{definition}
The above definition can be extended to a finite set of smooth maps $f_i:P_i\to M$ from manifolds with corners $P_i$ to a smooth manifold $M.$
\vskip .1in

\begin {theorem}\label{TH3.2}
1.  If  $f \pitchfork  S$    then $f^{-1}(S)$ is a smooth manifold  with corners (smooth submanifold  of $P$) with  
 $\partial _k f^{-1}(S)    =(f |_{\partial_kP} )^{-1}(S).$ Moreover if $P$ is clean then so is $f^{-1}(S).$   

2. If $f_1 \pitchfork f_2$  then $$(f_1\times f_2)^{-1}(\Delta_M):=  \{(u,v)\in P_1\times P_2 | f_1(u)= f_2(v)\}$$
 is a smooth manifold with corners. Moreover if $P_i$ are clean so is  $(f_1\times f_2)^{-1}(\Delta_M).$

\end{theorem}

Suppose $P$ is an oriented clean smooth manifold with corners.  An orientation of $P$ induces an orientation on $\partial _1 P$ and therefore an orientation on each component of $\partial _1 P.$ Let $\alpha_1$ and $\alpha_2$ be two components of $\partial_1 P$ and $\beta$ a component of $\partial_2 P$. Suppose that  and that 
$\beta\subset \overline \alpha_1$ and $\beta\subset \overline \alpha_2.$ Clearly $\beta$ is a codimension one submanifold of the   $\partial_1 P.$ Then an orientation $ o $ on $P$ induces an orientation $o_{\alpha_1}$ on $\alpha_1$  which in turn induces an orientation $o^1_\beta $ on $\beta.$ Similarly the orientation on $P$ induces an orientation $o_{\alpha_2}$ on $\alpha_2$ which in turn induces an orientation $o^2_\beta$ on $\beta.$
Observe that the orientations  $o^1_\beta$ and $\o^2_\beta$ are opposite. 

Suppose now that $P$ is a compact orientable clean smooth manifold with corners  of dimension $n.$

Fix an orientation $o_\alpha$ on each component $\alpha$ of  the corners.  Each  component $\alpha$ of $\partial_kP$ is  a smooth manifold of dimension $n- k.$ Denote by $\mathcal P_k$ the set of components of dimension $k$ (the  components of $\partial _{n-k} P$). 

Let $I_k :\mathcal P_k\times \mathcal P_{k-1}\to \mathbb Z$ be the  function defined by 

 \begin{equation*}
I_k(\alpha, \beta)=      { \begin{cases} 
0  \ if  \beta\nsubseteq \overline \alpha&\\
+1\  if \beta\subset \overline \alpha&\text{ and}\   o_\alpha\ \text{ induces} \ o_\beta\\
-1 \  if \beta\subset \overline \alpha&\text{ and}\   o_\alpha \ \text{ induces the oposite of } o_\beta
 \end{cases}}
 \end{equation*}
  $\alpha \in \mathcal P_k, \beta\in \mathcal P_{k-1}.$
 
The above observation implies that 
\begin{equation}\label {E;5}
\sum_{\beta\in \mathcal P_k} I_{k+1}(\alpha', \beta)\cdot I_k(\beta, \alpha'')=0
\end{equation}
for any $\alpha' \in \mathcal P_{k+1},$
$\alpha'' \in \mathcal P_{k-1}.$

If for a commutative ring $\kappa$ one considers the $\kappa-$module $$C^k: = \text{Maps}(\mathcal P_k, \kappa)$$ and 
the linear maps 
$d^k: C^k\to C^{k+1}$ defined by 
$$d^k(f)(\alpha)= \sum_{\beta\in \mathcal P_k}  I_{k+1}(\alpha, \beta)  f(\beta)$$
$f\in C^k,$ then the equality (\ref{E;5}) implies that $(C^k, d^k)$ is a cochain complex.

The cohomology of this cochain complex is called the {\bf incidence cohomology} of  the manifold with corners $P.$ This  cohomology is  independent on the chosen orientations $o_\alpha'$s.

\section{Proof of the main theorem}

We will prove Theorem \ref{TH1.5}  only for $\hat W^-_x$ and $\hat {\mathcal T}(x,y) .$ The statements  for $\hat W^+_x$ will follow by changing $X$ into $-X$ since the stable sets for $X$ are the unstable sets for $-X.$  The statements  for  $\hat {\mathcal M}(x,y)$ can be verified in essentially the same way as for  $\hat {\mathcal T}(x,y).$
Alternatively they  can be  derived  from the statements about $\hat W^-_x $ and 
$\hat W^+_y$ in view of the fact that $\hat{\mathcal M}(x,y)= (\hat i^-_x \times \hat i^+_y)^{-1}(\Delta_M)$  consists of the pairs of points in the product $\hat W^-_x \times \hat W^+_y$ equalized by $\hat i^-_x$ and $\hat i^+_y.$ 

We will focus the attention to the assertions (1.) and (2.).  Part of the assertion (3.), the orientability and the existence of the stable framing  follow from the simple observation that $\mathcal T(x,y)\times W^-_y$ is  an open set of $\partial_1 \hat W^-_x .$ The last part, the compatibility of the stable framings and of the orientability of $\mathcal T(x,y)$ for various $x,y$ is a tedious but conceptually straightforward verification. More details will be provided in the expanded version of this work,  cf. \cite {BH03},  which treats the more general case of Bott - hyperbolic vector fields.

We first prove assertions (1.) and (2.)
under the additional hypothesis  {\bf H}.

{\bf Hypothesis H:}  $X$  admits a proper Lyapunov function $f:M\to \mathbb R$ which in the neighborhood of rest points, in convenient coordinates $u_1,\cdots, u_k, v_1,\cdots, v_{n-k}$, is given by the quadratic expression

\begin {equation*}
f(x_1, \cdots x_n)= c -\sum_{i=1}^k u_i^2 + \sum_{j=1}^{n-k} v_j^2.
\end{equation*}
 We also suppose that with respect to these coordinates  
the unstable resp. the stable set of the rest points corresponds to $\mathbb R^k\times 0$ resp. $0\times \mathbb R^{n-k}.$ 
\vskip .1in

Let  $\cdots c_1 <c_2 < \cdots c_k< c_{k+1}\cdots $ be the set of  critical values.  Choose $\epsilon_i$ such that 
\ \ \ $c_{i+1}-\epsilon _{i+1} > c_{i}+\epsilon _{i}.$ 

Introduce 
\begin{equation*}
 M_i= f^{-1}(c_i), \ \   M_i^{\pm}= f^{-1}(c_i \pm \epsilon_i), \ \ M(i)= f^{-1}(c_{i-1}, c_{i+1}) 
 \end{equation*}
 and denote by $\mathcal X(i)$ the set of rest points which lie in $M_i,$  $ \mathcal X(i)= \mathcal X\cap M_i.$
Denote by 
\begin{equation*}
\varphi_k : M_{k}^-\to M^+_{k-1}
\end{equation*}
the map defined by the flow of $X.$ Precisely $\varphi_k(x)$ is the intersection of the trajectory $\gamma_x,\ x\in M^-_k$  with  $M^+_{k-1},$
see FIGURE 2. below.
\begin{figure}[ht]
\begin{center}
\includegraphics[height=90mm]{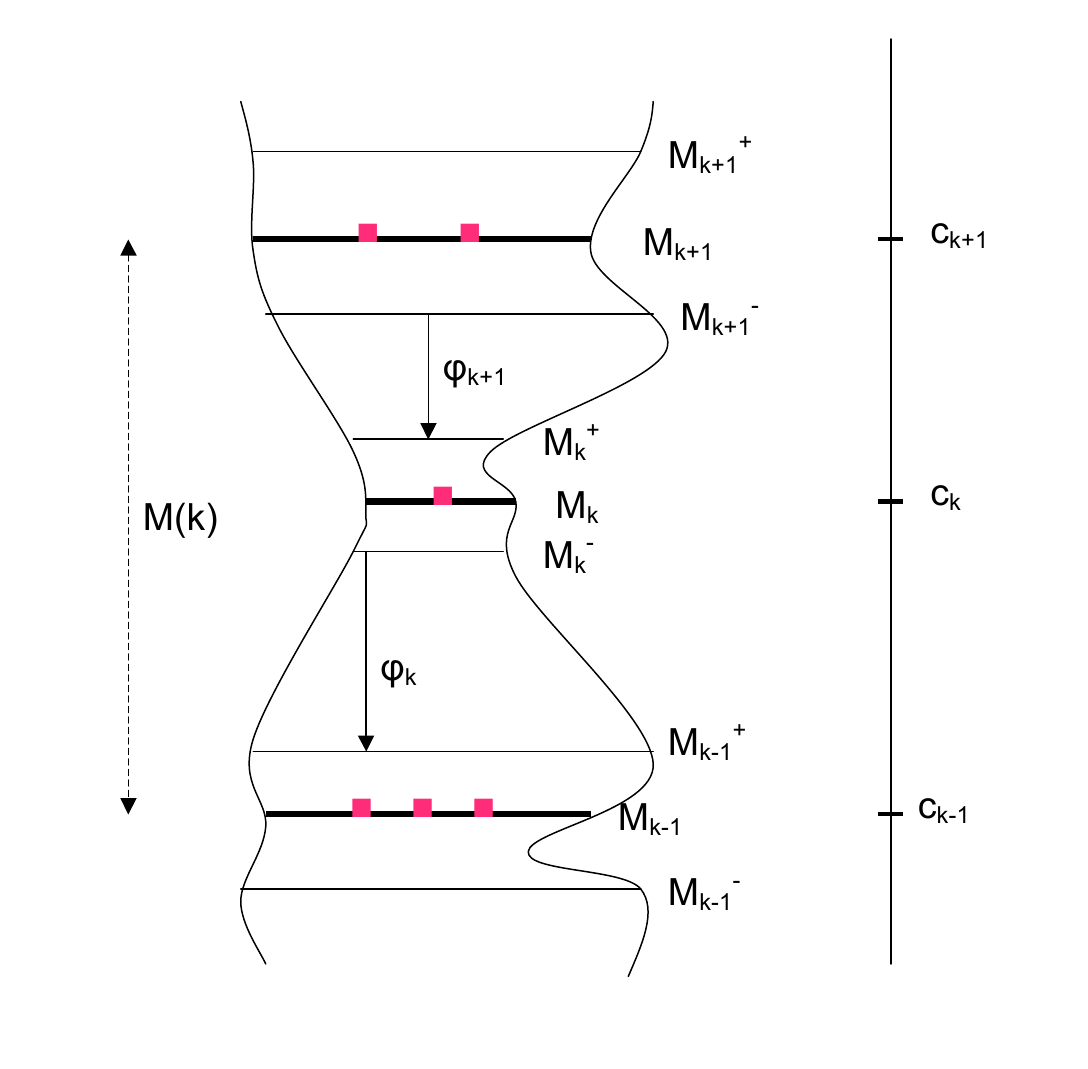}
\caption{}
\label{fig:DANU7a}
\end{center}
\end{figure}

For $x\in \mathcal X(i)$ denote by 
\[S^\pm_x= M^\pm_i\cap W^\pm_x, \ \ D^\pm_x= M(i)\cap W^\pm_x,\]
see FIGURE 3. below,
\[\mathbb S^{\pm}_i:= \bigcup _{x\in\mathcal X(i)} S^\pm_{x},\]

\[\mathbb D^{\pm}_i:= \bigcup _{x\in\mathcal X(i)} D^\pm_{x},\]

\[(\mathbb S^+\rtimes \mathbb S^-)_i:= \bigcup _{x\in\mathcal X(i)} S^+_{x}\times S^-_x.\]

\[(\mathbb S^+\rtimes \mathbb D^-)_i:= \bigcup _{x\in\mathcal X(i)} S^+_{x}\times D^-_x.\]
Introduce the sets 
\[P_i=\{(u,v) \in M^+_i\times M^-_i  \mid 
 u,v \ \text{lie on the same possibly broken  trajectory} \},\] 
\[Q_i=\{ (u,v ) \in M^+_i\times M(i) \mid 
u, v \  \text{ lie on the same possibly broken  trajectory}\}.\] 

\begin{figure*}[ht]
\begin{center}
\includegraphics[height=50mm]{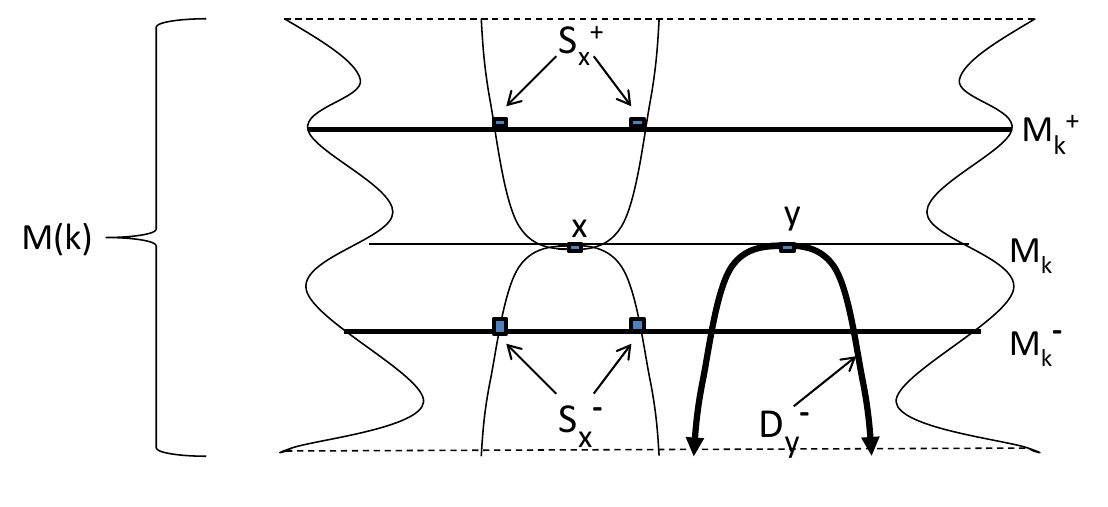}
\caption{}
\label{fig:DANU7a}
\end{center}
\end{figure*}

Each  set  $P_i$ resp. $Q_i$ is a union of two disjoint subsets 
$\overset {\circ}{P_i}$ and $\partial P_i$ resp. $\overset {\circ}{Q_i}$ and $\partial Q_i,$ 
\[\overset{\circ} P_i= \{(u,v) \in P_i \mid u,v \ \text{ lie on the same  unbroken trajectory} \} \]
\[\partial P_i= \{ (u,v) \in P_i \mid 
u,v \ \text{ lie on the same  broken trajectory}\}\]
\[\overset{\circ} Q_i= \{(u, v) \in P_i \mid u,v \ \text{lie on the same  unbroken trajectory}\}\]
\[\partial Q_i= \{(u,v)  \in Q_i \mid u,v \ \text{ lie on the same  broken trajectory}\}.\]

Equipped  with the  topology induced from the product $M^+_i\times M^-_i$ resp. $M^+_i\times M(i)$
 $P_i$ resp. $Q_i$ is a  topological manifold with boundary
with  $\overset{\circ} P_i \ \text {resp.}\  \overset{\circ} Q_i$  the interior and $\partial P_i \ \text {resp.}\  \partial  Q_i$
 the boundary.  Actually both   $\overset{\circ} P_i\  \text{and}\  \partial P_i $ resp.$\ \overset{\circ} Q_i\  \text {and}\  \partial  Q_i$
 are smooth submanifolds of $M^+_i\times M^-_i$ resp.  $M^+_i\times M(i)$  but $P_i$ resp. $Q_i$ might not be  in general.
 
For both $P_i$ and $Q_i$ denote by $p^+_i$ resp. $p^-_i$ the projection on the first resp. second component of  $M^+_i\times M^-_i$ and $M^+_i\times M(i)$ and  by $p_i= p^+_i\times p^-_i$ their product. The map $p_i$ is one to one and homeomorphism onto  the image. In addition we have the following.
\begin{enumerate} 
\item The map $p_{i}$ identifies $\partial P_i$ with $(\mathbb S^+\rtimes \mathbb S^-)_i,$ and  $\partial Q_i$ with $(\mathbb S^+\rtimes \mathbb D^-)_i.$  
\item The  restrictions of $p^\pm_{i}$ to $\overset {\circ}{P_i}$  are  diffeomorphisms 
 onto $M^\pm_i\setminus \mathbb S^\pm_i$ and the restrictions to $\partial P_i,$  via the identification above,  
  are the projections on  $\mathbb S^\pm_i.$ 
 \item The restriction  of $p^+_{i}$ to $\overset {\circ}{Q_i}$  is a smooth bundle over $M^+_i\setminus \mathbb S^+_i$ with fiber an open interval  and the restriction to $\partial Q_i$  
 is the projection on $\mathbb S^+_i.$
\item The restriction of $p^-_{i}$ to $\overset {\circ}{Q_i}$ is a diffeomorphism onto $M(i)\setminus \mathbb D^-_i$ and the restriction to $\partial Q_i$  
is the projection on $\mathbb D^-_i.$ 
\end{enumerate} 
As pointed out  above the subsets 
$p_{i}(P_i)$ and $p_{i}(Q_i)$ are not necessarily smooth submanifolds of $M^+_i\times M^-_i$ resp.
$M^+_i\times M(i),$  
however the  following propositions will provide  structures of smooth   manifold with boundary on both $P_i$ and  $Q_i$ which will make $p_{i}$  smooth maps.

\begin{proposition}\label {P:a}
The smooth map   $p_{i}: \partial P_i\to   M^+_i\times M^-_i$ admits  smooth extensions $ \tilde p_{i} :\partial P_i\times [0,\epsilon)\to M^+_i\times M^-_i$ so that:

1. the image is  a neighborhood of $(\mathbb S^+\rtimes \mathbb S^-)_i $ in $p_{i}(P_i),$

2. $\tilde p_{i}$ is injective,

3. $ \tilde p_{i}$ restricted to $\partial P_i \times (0,\epsilon)$
is of maximal rank.
\end{proposition}
\begin{proposition} \label {P:b}
The smooth map   $p_{i}: \partial Q_i \to   M^+_i\times M(i)$ admits  smooth extensions $\tilde p _{i} :\partial Q_i\times [0,\epsilon)\to M^+_i\times M(i)$ so that:

1. the image is  a neighborhood of $(\mathbb S^+\rtimes \mathbb D^-)_i $ in $p_{i}(Q_i)$

2. $\tilde p_{i}$ is injective.

3. $\tilde p_{i}$ restricted to $\partial Q_i \times (0,\epsilon)$
is of maximal rank.
\end{proposition}
The proofs of these propositions will be given towards the end of the section.

Equip $P_i$  with the smooth structure defined by the atlas obtained from  $\{\partial P_i\times [0,\epsilon), \tilde p_{i}\}$ 
 and 
 $\overset{\circ}{P_i}.$ Similarly  equip  $Q_i$ with the smooth structure defined by the atlas obtained from  $\{\partial Q_i\times [0,\epsilon), \tilde p_i\}$ 
 and 
 $\overset{\circ}{Q_i}.$    Equivalently,  regard $P_i$ resp. $Q_i$ obtained by glueing $\partial P_i\times [0,\epsilon)$ resp. $\partial Q_i\times [0,\epsilon)$ to $\overset{\circ}{P_i}$ resp. $\overset{\circ}{Q_i}$ via the diffeomorphisms provided by the restriction of $\tilde p_{i}$  to $\partial P_i\times (0,\epsilon)$ resp. to $\partial Q_i\times (0,\epsilon).$ These smooth structures will be denoted by $(P_i)_h$ resp.$(Q_i)_h.$ 
 
If the rest points of $X$ are of Morse type  then we have the following.
\begin{proposition} \label{P:c}  
 If the rest points of $X$ are of Morse type  then 
 the image $ p_{i}(P_i)\subset M^+_i\times M^-_i$ resp. $p_{i}(Q_i)\subset M^+_i\times M(i)$ are  smooth submanifolds with boundaries.
\end{proposition}  
The proof of this proposition  will be given towards the end of the section.
\vskip .1in

This implies that 
$P_i$ resp. $Q_i$ have a smooth structure of manifold with boundary  denoted by $(P_i)_m$ resp. $(Q_i)_m.$ The structures $(m)$  are  $(h)$ are never   the same but  $id: (P_i)_h\to (P_i)_m$ and  $id: (Q_i)_h\to (Q_i)_m$
are smooth homeomorphisms which restrict to diffeomorphisms on  the interiors and on  the boundaries.
\vskip .1in

Propositions  \ref {P:a},  \ref {P:b} and \ref {P:c} 
imply   that for any $(r,k)$ the  product  $\mathcal P:= A\times P_{r+k-1}\times P_{r+k-2}\cdots P_{k+1}\times B$ with $A$ a smooth manifold  and $B$ a smooth manifold, possibly  with boundary,
is a smooth manifold with corners.  The   corner   
$\partial _l \mathcal P$ can be described as follows.

For any $i$ with $r+k-1 \geq i\geq k+1$ denote by $R_i$  the subset of $P_i$ which is either the interior or the boundary of  $P_i$  and by 
$R_k$ the subset of $B$ which is either the interior or the boundary  of $B.$ 

Then the corner $\partial _l \mathcal P$ is the disjoint union of products $ A\times R_{r+k-1}\times R_{r+k-2}\cdots R_{k+1}\times R_{k}$ with $l$  of the sets $R_i's$ being  boundaries  
and the remaining $r-l$   being  interiors.  For example  if $B$ is a smooth manifold with boundary  and $l=1$ 
\begin{equation*}
\begin {aligned}
\partial_1 \mathcal P=& A\times \partial P_{r+k-1}\times \overset{\circ} P_{r+k-2}\times \cdots \overset{\circ} P_{k+1}\times \overset {\circ}B \sqcup \\
&A\times \overset{\circ} P_{r+k-1}\times  \partial  P_{r+k-2}\times \cdots \overset{\circ} P_{k+1}\cdots  \times \overset {\circ}B \sqcup \cdots \\
&A\times \overset{\circ} P_{r+k-1}\times \overset{\circ} P_{r+k-2}\times \cdots \partial P_{k+1}\times \overset {\circ}B\sqcup\\
&A\times \overset{\circ} P_{r+k-1}\times \overset{\circ} P_{r+k-2}\times \cdots \overset{\circ} P_{k+1}\times \partial B.
\end{aligned}
\end{equation*}
Suppose Propositions \ref{P:a}, \ref{P:b}, \ref{P:c}  were established. Here is the general scheme to  verify 
(1.) and (2.) 
for $\hat{\mathcal T}(x,y)$ and $\hat W^-_x.$  

For any $r,k\geq 0$  consider the diagram 

\xymatrix{
X_{r+k} \ar[d]_{id}\ar[dr]^{\varphi_{r+k-1}}& & X_{r+k-1} \ar[d]_{id}\ar[dr]^{\varphi_{r+k-2}} & & X_{r+k-2} \ar[d]_{id} \cdots  & & &\\
X_{r+k}&Y_{r+k-1}&X_{r+k-1}&Y_{r+k-2}& X_{r+k-2}\cdots & & &
\\
A \ar[u] _\alpha& & Z_{r+k-1}\ar[ul]^{p^+_{r+k-1}} \ar[u]_{p^-_{r+k-1}}& & Z_{r+k-2}\ar[ul]^{p^+_{r+k-2}} \ar[u]_{p^-_{r+k-2}} \cdots  & & & 
}
\vskip .1in

\xymatrix{
\cdots  &&  X_t \ar[d]_{id}\ar[dr]^{\varphi_{t-1}}& & X_{t-1}\ar[d]_{id}\ar[dr]^{\varphi_{t}} &\ \ &\cdots  X_{k+1} \ar[d]_{id}\ar[dr] ^{\varphi_{k}}&&\\
\cdots  &&  X_t& Y_{t-1}&X_{t-1} &\ \ &\cdots  X_{k+1} & Y_k\\
 \cdots &&   Z_t\ar[ul]\ar[u]_{p^-_t}& & Z_{t-1} \ar[ul]^{p^+_{t-1}} \ar[u]_{p^-_{t-1}} &\ \ &\cdots Z_{k+1}\ar[u]_{p^-_{k+1}} & B\ar[u]_\beta
}
\vskip .2in
\hskip 2in {FIGURE 4.}
\vskip .2in

\noindent where $X's, Y's$ and $A$ are  smooth  manifolds, $Z's$ and $B$   smooth manifolds  with boundary  (possibly empty) and the  arrows  are smooth maps with $\varphi'$s embeddings. 

Denote by  $\mathcal S, \mathcal O$ and $ \mathcal P$ the spaces defined by:
\begin{enumerate}
\item $\mathcal S= \mathcal S(r+k,k) := X_{r+k}\times X_{r+k-1} \times \cdots X_{k+1},$
\item
$\mathcal O= \mathcal O(r+k,k):= X_{r+k}\times Y_{r+k-1}\times X_{r+k-1}\times \cdots X_{k+1}\times Y_k,$
\item
$\mathcal P=\mathcal P(r+k,k):= A\times Z_{r+k-1}\times Z_{r+k-2} \times \cdots Z_{k+1}\times B,$
\end{enumerate}
Denote by $s$ and $t$  the maps defined by: 
\begin{enumerate}
\item
$s= s(r+k,k): \mathcal S\to \mathcal O,$ 
the product of all maps from the top line to the  middle line in 
the diagram above and 
\item
$t=t(r+k,k): \mathcal P\to \mathcal O,$ 
 the  product of all maps from the bottom  to the middle line in diagram above.
\end{enumerate}
Since $id$ and $\varphi's$ are embeddings so is $s$ and $s(\mathcal S)$ is a smooth submanifold of $\mathcal O.$ 
Note that   $\mathcal P$ is a smooth manifold with corners, $\mathcal O$ and $ \mathcal S$ are smooth manifolds,  $t$ and $s$  are smooth maps. 

We  say  that " the diagram FIGURE 4. is transversal"  if 
$t \pitchfork s(\mathcal S).$   If so   by Theorem \ref {TH3.2},  the subspace  $t^{-1}(s(\mathcal S))$ receives  a structure of  smooth manifold with corners.

To accomplish the proof of (1.) and (2.) in Theorem \ref{TH1.5} we  choose the integers $r, k,$   
the manifolds  $X's, \cdots, B$ and  the maps $\varphi_i, p^\pm_i, \alpha, \beta, $ appropriately in order to  obtain  $\hat {\mathcal T}(x,y)$  and  
 $(\hat i_x)^{-1}(f^{-1}(c_{r+k}, c_r)),$ 
as   $t^{-1}(s(\mathcal S)).$ Then we verify the transversality of the diagram FIGURE 4.

{\it  The case  of $\hat{\mathcal T}(x, y).$ } 
Choose $(r,k)$ so that $f(x)= c_{r+k}$ and $f(y)=c_k.$

Take   
$A= S^-_x,$ $B= S^+_y,$
\begin{equation*}
\begin{aligned}
(\varphi_{i+1}: X_{i+1}\to Y_{i})=&(\varphi_{i+1}:  M^-_{i+1} \to M^+_{i}), \   \    k \leq i\leq r+k-1,  \\
(Z_j, p^+_j, p^-_j)=&(P_j, p^+_j, p^-_j), \  \ k+1\leq j \leq r+k-1.
\end{aligned}
\end{equation*}
Take 
$\alpha, \beta$ to be the obvious inclusions.  
 With these choices  
 $f^{-1}(s(\mathcal S))$ identifies with $\hat {\mathcal T}(x,y).$  
 
Diagram FIGURE 4 becomes 

\xymatrix{
M^-_{r+k} \ar[d]_{id}\ar[dr]^{\varphi_{r+k-1}}& & M^-_{r+k-1} \ar[d]_{id}\ar[dr]^{\varphi_{r+k-2}} & & M^-_{r+k-2} \ar[d]_{id} \cdots  & & &\\
M^-_{r+k}&M^+_{r+k-1}&M^-_{r+k-1}&M^+_{r+k-2}& M^-_{r+k-2}\cdots & & &
\\
{S^-_x}\ar[u] _\alpha& &P_{r+k-1}\ar[ul]^{p^+_{r+k-1}} \ar[u]_{p^-_{r+k-1}}& & P_{r+k-2}\ar[ul]^{p^+_{r+k-2}} \ar[u]_{p^-_{r+k-2}} \cdots  & & & 
}
\vskip .1in

\xymatrix{
\cdots &&M^-_t \ar[d]_{id}\ar[dr]^{\varphi_{t-1}}& & M^-_{t-1}\ar[d]_{id}\ar[dr]^{\varphi_{t}} &\ \ &\cdots  M^-_{k+1} \ar[d]_{id}\ar[dr] ^{\varphi_{k}}&&\\
\cdots &&M^-_t& M^+_{t-1}&M^-_{t-1} &\ \ &\cdots  M^-_{k+1} &M^+_k\\
\cdots &&P_t\ar[ul]\ar[u]_{p^-_t}& & P_{t-1} \ar[ul]^{p^+_{t-1}} \ar[u]_{p^-_{t-1}} &\ \ &\cdots P_{k+1}\ar[u]_{p^-_{k+1}} & S^+_y\ar[u]_\beta
}

\vskip .1in
\hskip 2in {FIGURE 5.}
\vskip .2in
      
Then  (1.) and  (2.)  follow  from Proposition \ref{P:d} below.
 
 \begin{proposition}\label {P:d}
 The diagram (FIGURE 5) is transversal. 
\end{proposition} 
Propositions \ref {P:d} is a consequence of the transversality $i^-_x\pitchfork i^+_y$ for $x,y\in \mathcal X(X).$ For details the reader can consult 
\cite {BFK} and \cite {BH01}.

{\it The case of  $\hat {W}^-_x:$}  Let $f(x)= c_m$ and  $k<m.$ 
This case is treated in two steps. 

First we check  that 
 the open set $(\hat i_x)^{-1}(M(k))$ 
 has a structure of a smooth manifold with corners.  
For this purpose we use the same diagram (FIGURE 4) for $(r=m-k, k)$ $X, Y, Z, A,\alpha $  as in the case of  $\hat{\mathcal T}(x, y),$ $B= Q_k,$   $\beta = p^+_{k}: Q_k \to M^+_k$ and we replace Proposition 
\ref{P:d} by  Proposition \ref {P:e} below.
Diagram  FIGURE 4 becomes 

\xymatrix{
M^-_{r+k} \ar[d]_{id}\ar[dr]^{\varphi_{r+k-1}}& & M^-_{r+k-1} \ar[d]_{id}\ar[dr]^{\varphi_{r+k-2}} & & M^-_{r+k-2} \ar[d]_{id} \cdots  & & &\\
M^-_{r+k}&M^+_{r+k-1}&M^-_{r+k-1}&M^+_{r+k-2}& M^-_{r+k-2}\cdots & & &
\\
{S^-_x}\ar[u] _\alpha& & P_{r+k-1}\ar[ul]^{p^+_{r+k-1}} \ar[u]_{p^-_{r+k-1}}& & P_{r+k-2}\ar[ul]^{p^+_{r+k-2}} \ar[u]_{p^-_{r+k-2}} \cdots  & & & 
}
\vskip .1in

\xymatrix{
&&M^-_t \ar[d]_{id}\ar[dr]^{\varphi_{t-1}}& & M^-_{t-1}\ar[d]_{id}\ar[dr]^{\varphi_{t}} &\ \ &\cdots  M^-_{k+1} \ar[d]_{id}\ar[dr] ^{\varphi_{k}}&&\\
&&M^-_t& M^+_{t-1}&M^-_{t-1} &\ \ &\cdots  M^-_{k+1} & M^+_k\\
&&P_t\ar[ul]\ar[u]_{p^-_t}& & P_{t-1} \ar[ul]^{p^+_{t-1}} \ar[u]_{p^-_{t-1}} &\ \ &\cdots P_{k+1}\ar[u]_{p^-_{k+1}} & Q_k\ar[u]_{p^+_k}
}
\vskip .1in
\hskip 2in {FIGURE 6.}
\vskip .2in

\begin{proposition} \label {P:e}
The diagram (FIGURE 6) is transversal.
\end{proposition}
As with  Proposition \ref {P:d},  Proposition \ref{P:e} is a  consequences of the transversality $i^-_x\pitchfork i^+_y$ for $x,y\in \mathcal X(X).$ For details the reader can consult 
\cite {BFK} and \cite {BH01}.

Proposition \ref {P:e}  implies  that $(\hat i_x)^{-1}(M(k)) $ has a structure of smooth manifold with corners. 

Second, we verify  that the smooth structures on $(\hat i_x)^{-1}(M(k)) $ and on $(\hat i_x)^{-1}(M(k-1)) $ agree.
For this purpose we 
consider the  map  $h:= f\circ p^-_i: Q_i\to (c_{i+1},  c_{i-1})$ 
and let 
$Q_i':=  h^{-1}(c_{i+1}, c_i)$ and  $  Q_i'':= h^{-1}(c_i, c_{i-1}).$ Both are open subsets of $Q_i$
and we have:
\begin{observation}\label {O:a}
There are canonical diffeomorphisms 
$$\theta'_k: Q_k'\to M^+_k\times (c_{k+1}, c_k)$$
$$\theta''_k: Q_k'' \to P_k\times_{M^-_k} (M^-_k\times (c_k, c_{k-1})) $$ where the fiber product is taken with respect to $p^-_{k}: P_k\to M^-_k$ and the projection $M^-_k\times (c_k, c_{k-1})\to M^-_k.$
\end{observation}

Then  the composition  of :
\begin{enumerate}
\item $\theta_k'' :Q_k'' \to P_k\times _{M^-_k} (M^-_k\times (c_k, c_{k-1})),$
\item the inclusion $ P_k\times _{M^-_k} (M^-_k\times (c_k, c_{k-1}))\subset P_k\times  (M^-_k\times (c_k, c_{k-1})),$
\item $ id\times \varphi_{k-1}\times id : P_k\times  M^-_k\times (c_k, c_{k-1}) \to P_k\times  M^+_{k-1}\times (c_k, c_{k-1})$ and 
\item $id\times (\theta_{k-1}')^{-1}: P_k\times  M^+_{k-1}\times (c_k, c_{k-1})\to P_k\times Q_{k-1}'$ 
\end{enumerate}
\noindent is a smooth embedding denoted by 
$\theta_k: Q_k'' \to P_k\times Q_{k-1}'.$ 

We write $t':\mathcal P'\to \mathcal O$ resp. $t'': \mathcal P'' \to \mathcal O$ instead of  $t:\mathcal P\to \mathcal O$ and  $Q'$ resp. $Q''$ instead of $Q.$ 

In view of Observation \ref{O:a}  the map 
$  \mathcal P'(n,k)\to \mathcal P''(n,k-1),$ given by the product of $id's$ (on $ S^-_{x}, P_{n-1},\cdots P_{k-1}$)  and of $\theta_k: Q_k'' \to P_k\times Q_{k-1}',$ 
 is a smooth embedding which sends 
$$(t'' (n,k))^{-1}(s(n,k)(\mathcal S(n,k))$$ onto $$(t' (n,k-1))^{-1}(s(n,k-1)(\mathcal S(n,k-1)) .$$  It identifies  the structures of smooth manifolds with corners  which were derived using $k$ and $k-1.$ 

Apparently the smooth structures defined so far depend on the Lyapunov function and the choices of $\epsilon_i;$ this is not the case.

\vskip .1in
{\bf The independence of Lyapunov function:} The arguments  are the same for $\hat {\mathcal T}$ and $\hat{W}^-$ so we will treat  only $\hat {\mathcal T}.$

If $f$ and $f'$ are two Lyapunov functions and $\underline \gamma$  a  possibly broken instanton  we consider   two transversals $\underline V= (V_1,\cdots, V_k)$  and $\underline V'=(V'_1,\cdots V'_k)$ with the same mark points  and  $V_i$ contained in the levels of $f$ and $V'_i$ contained in the levels of $f'.$ We can find a diffeomorphism  $\theta$ of a neighborhood $U$ of $\gamma$ onto a  neighborhood $\theta(U)$ of $\gamma$ which restricts to the identity on $\gamma$  and sends $V_i$ into $V'_i.$  This diffeomorphism provides an open embedding from the product of $V_i$ into the product of $V'_i.$ 
It follows that $id: \hat {\mathcal T}\to \hat {\mathcal T}$ is smooth and of maximal rank in the neighborhood of $\underline \gamma$
with respect to either one of the smooth structure  $(h)$ or $(m)$ defined using $f$ and $f'.$ 
\vskip .1in 
{\bf The  removal of the additional hypothesis} {\bf H:}  While global Lyapunov functions might not exist,   for any broken instanton  $\underline \gamma$ from the rest point $x$ to the rest point $y$ 
one can find an open neighborhood $U$ in $M$  so that  a "convenient Lyapunov function"  $f: U\to \mathbb R$ for $X|_U$ exists.  
Here "convenient  Lyapunov function"  means that the system $(X|_U, f, U)$ is diffeomorphic to $(Y|_V, g|_V, V)$ where $Y$ is a smooth vector field satisfying $P_1$ and $P_2$ on a smooth  manifold $N,$  $g:N\to \mathbb R$ a proper Lyapunov function for $Y$ and $V$ an open set in $N.$
As  the space $\hat{\mathcal T}_U(x,y)$ consisting of broken instantons from   $x$ to $y$ which lie in $U$ is an open set in $\hat{\mathcal T}(x,y)$  we define a smooth structure on  $\hat{\mathcal T}_U(x,y)$ and 
note  that for different such $U's$ these structures  agree on intersections. 
The smooth structure on $\hat{\mathcal T}_U(x,y)$ is defined using the space of broken instantons of $Y$ on $N$ which lie in $V.$
\vskip .1in

{\it Proof of  Proposition \ref{P:a}, \ref {P:b} }

First we introduce some notation. In the context  of Theorem  \ref {T;2.1} in section 2
denote by $ S^\pm$ and by $ D^\pm$ the  sphere and  the disc of  radius $\epsilon$ in  $\mathbb R^{\pm}$ and when this notation is applied to coordinates about a rest point $x$ write $ S^\pm_x$ and by $ D^\pm_x$ instead. 

Define  the maps 
 $\chi_1= (\chi^+_1, \chi^-_1):  S^+\times  D^-\times [0, \infty )\to (\mathbb R^+\times \mathbb R^-)$ 
and 
 $\chi_2= (\chi^+_2, \chi^-_2):  S^+\times  D^-\times [0, \infty)\to (\mathbb R^+\times \mathbb R^-),$ 
by  the formulae:
\begin{equation}\label {F;c}
\begin{aligned}
\chi_1^+(p, q, s)=&p\\
\chi_1^-(p, q, s)=&\gamma^-(-1/s, p, q, -1/s, 1/s )\text \ if s\ne 0\\
\chi_1^-(p, q, 0)=&0
\end{aligned}
\end{equation}
\begin{equation}\label {F;d}
\begin{aligned}
\chi_2^+(p, q, s)=&\gamma^+(1/s, p, q, -1/s, 1/s
)\text \ if s\ne 0\\
\chi_2^-(p, q, 0)=&0\\
\chi_2^-(p, q, s)=&q. 
\end{aligned}
\end{equation}

Define  $\chi:=(\chi_1, \chi_2):  S^+\times  D^-\times [0, \infty)\to (\mathbb R^+\times \mathbb R^-)\times  
(\mathbb R^+\times \mathbb R^-).$

Clearly $\chi ( S^+\times S^-\times[0,\infty))\subset ( S^+\times \mathbb R^-)\times (\mathbb R^+\times \mathbb S^-).$ 

The estimates in Theorem \ref {T;2.1} 
show  that the map $\chi$ is smooth  and 
for $\theta $  small  the restriction of  $\chi$  to $S^+\times S^-\times(0,\theta)$ 
and  to
$ S^+\times S^-\times 0$ is of maximal rank but $\chi$  is not.
It fails  at the points of $ S^+\times  S^-\times 0
.$ 
\vskip .1in 
{\bf Example :} If 
\begin{equation}\label{Eq;1}
\begin{aligned}
X =- \grad f= -&\sum_{i=1}^k  x_i\partial / \partial x_j +  \sum_{j=k+1}^n
x_j\partial / \partial x_j.\\
p= (x_1,\cdots x_k), \  & \  q=(x_{k+1}, \cdots, x_n)
\end{aligned}
\end{equation} 
a simple calculation shows that:
$$\gamma^+(t;  p, q, T_1,T_2)= e^{(T_1-t)}p, \ \   \gamma^-(t;  p, q, 0,T)= e^{-(T_2-t)}q .$$
The  estimates in Theorem \ref {T;2.1}   are  satisfied and  $\chi$ is visibly  not of maximal rank at the points of $ S^+\times  S^- \times 0.$ 
\vskip .1in

We proceed now with the proof of Propositions \ref{P:a} and \ref {P:b}.

Observe that it suffices to check the statements in Propositions \ref{P:a} and \ref {P:b} for $\epsilon _i$ small enough and  the statement in Proposition  \ref {P:b}  for  $M(i)$ replaced by the smaller open set $f^{-1}(c_i-\epsilon_i, c_i + \epsilon_i).$

Choose for each rest point $x\in\mathcal X(i)$ a  neighborhood and coordinates in the neighborhood   so that the hypotheses of Theorem \ref {T;2.1}
are satisfied and  $f = c_i -1/2 |p|^2 + 1/2|q|^2.$ Here  $|p|$ and $|q|$ denote the norm in the respective  coordinates.
Choose $\epsilon_i =\epsilon/2 $  with $\epsilon$ small enough to  have the conclusions of Theorem  \ref {T;2.1} satisfied for  each rest point. 

Since there is no risk of confusion from now on we drop the index $i$ from notation and write $M(\epsilon)$ instead of $f^{-1}(c_i-\epsilon, c_i + \epsilon).$

Define 
  $u^+_x:  S^+_x\times  D^-_x\to M^+, \text{resp}.\   
u^-_x:  D^+_x\times  S^-_x \to M^-$ to be the map which assigns to $(p,q)$ the intersection of the trajectory through $(p,q)$ with $M^+$ resp. $M^-.$  
The maps  $u^+_x$ and  $u^-_x$ are diffeomorphisms on their images
and their restrictions to $ S^+_x\times 0$ resp. to $0\times S^-_x$ are the identity maps. 

For the rest point $x$ denote by $(\chi^x_1, \chi^x_2)$ the maps $(\chi_1, \chi_2)$ defined by the formulae (\ref {F;c}) and (\ref  {F;d}).

For $P$ take  $\tilde p_x:=  (u^+_x\circ \chi^x_1, u^-_x\circ \chi^x_2):  S^+_x\times  S^-_x\times [0, \epsilon')\to M^+\times M^-$ and for $Q$ take $\tilde p_x:= (u^+_x\circ \chi^x_1,  \chi^x_2):  S^+\times  D^-_x\times [0, \theta)\to M^+\times M(\epsilon)$ with $\theta$ small enough to insure that the image of  $S^+_x\times  S^-_x\times [0, \theta)$ by $\chi^x_2$ lies in $ M(\epsilon).$   Take $\tilde p= \sqcup_{x\in \mathcal X(i)} \tilde p_x.$ The maps $\tilde p$ satisfy the conclusions of Propositions \ref {P:a} and \ref{P:b}. 
\vskip .1in

{\it Proof of Proposition  \ref{P:c}:} We use the same conventions and notations as in the previous proof.

The chosen  neighborhoods and coordinates
for  the rest points are so made  to have $X$  given by 
(\ref{Eq;1}). 
Then,  each 
 trajectory of  $X$   passing through $(p,q)\in \mathbb R^n\times \mathbb R^{n-k}$ ( $k= \ind x$)
 at $s=0$ is given by 
$$
\gamma(s)= (\gamma^+(s), \gamma^-(s)),  \ \  \gamma^+(s)= e^{-s} p,\ \ \gamma^-(s)= e^{s}q.$$

To check that $P$ is a smooth submanifold with boundary  it suffices to construct the  smooth maps 
$\omega_x: S^+ _x\times  S^-_x\times [0,\epsilon)\to M^+\times M^-$
so that $\omega : (\mathbb S^+\rtimes \mathbb S^-)\times [0,\epsilon) \to M^+\times M^-$  defined by $\omega= \sqcup_{x\in \mathcal X(i)} \omega_x$ satisfies:

1a.  $ \omega$ restricted to $\partial P= (\mathbb S^+ \rtimes \mathbb S^-)$ is the identity,

2a.  the image of $\omega$ is an open neighborhood of $\partial P$ in $P.$

3a.  $\omega$ is of maximal rank on $(\mathbb S^+ \rtimes \mathbb S^-) \times [0,\epsilon),$ 
(Note the distinction between item 3a. above  and item 3. in Propositions \ref{P:a} and \ref {P:b}.)

Define  $u_x= (u^+_x, u^-_x)$ by assigning to $(p,q,t),$ $p\in    S^+_x, q\in  S^-_x$ and $t\in [0,\epsilon)$  the pair of points provided by the intersection of the trajectory passing through $
p +tq,$  with $M^+$ and $M^-.$ Define $\omega_x:= (u^+_x\circ \chi^x_1, u^-_x\circ \chi^x_2)$ where  $$\chi^x(p,q,t):= (\chi^x_1(p,q,t), \chi^x_2(p,q,t):= (p, tq, tp, q).$$ Items  1a., 2a, 3a. above are  satisfied. 

 To check that $Q$ is a smooth submanifold with boundary  it suffices to construct the  smooth maps 
$\omega_x:  S^+_x \times  D^-_x\times [0,\epsilon)\to M^+\times M(\epsilon)$
so that $\omega= \sqcup_{x\in \mathcal X(i)}\omega_x$ satisfies:

1b.  $ \omega$ restricted to $\partial Q= (\mathbb S^+ \rtimes \mathbb D^-)$ is the identity,

2b.  the image of $\omega$ is an open neighborhood of $\partial Q$ in $Q,$

3b.  $\omega$ is of maximal rank on $(\mathbb S^+ \rtimes \mathbb D^-)\times [0,\epsilon).$ 

Define  $u^+_x$ by assigning to $(p,q,t),$ $p\in    S^+_x, q\in  D^-_x$  the intersection of the trajectory passing through $
tp +q$  with $M^+.$ 
Define  $\omega_x:= (u^+\circ \chi^x_1,  \chi^x_2)$ where  $$\chi(u,v,t)= (\chi^x_1(p,q,t), \chi^x_2(p,q,t)):= ((p, tq), (tp, q)).$$  Items 1b, 2b, 3b. above are  satisfied. 

q.e.d
 \begin{observation}
The reader can notice that the smooth structures (h) and (m) in the case of a vector field with the rest points of Morse type provided by Propositions \ref {P:a}, \ref {P:b} and Proposition \ref {P:c} respectively can not be the same. 
 \end{observation}

\end{document}